\documentclass[10pt]{amsart}
\usepackage{amssymb}
\usepackage{dsfont}
\usepackage{amscd}
\usepackage[mathscr]{euscript} 
\usepackage[all]{xy}
\usepackage{url}
\usepackage{stmaryrd}
\usepackage{comment}
\usepackage[retainorgcmds]{IEEEtrantools}
\usepackage{hyperref}
\usepackage{graphicx}
\usepackage{mathtools}
\usepackage{centernot}
\usepackage{marvosym}
\usepackage{tikz}
\usepackage{tkz-fct}
\usepackage{soul}
\usepackage{lineno}

\usepackage[T1]{fontenc}

\title[Of model completeness and algebraic groups] {Of model completeness and algebraic groups}

\author[D. M. HOFFMANN]{Daniel Max Hoffmann$^{\clubsuit}$}
\thanks{2020 \textit{Mathematics Subject Classification}.
Primary 03C10;
Secondary 03C60, 14L35}
\thanks{\textit{Key words and phrases}. semisimple algebraic groups, model completeness}
\thanks{$^{\clubsuit}$SDG. Supported by the European Union’s Horizon research and innovation
programme under the MSCA project no. 101063183,
by the National Science Centre (Narodowe Centrum Nauki, Poland)
grant no. 2021/43/B/ST1/00405, and by the UW IDUB PSP no. 501-D110-20-3004310.}
\address{$^{\clubsuit}$
Institut f\"{u}r Geometrie\\
Technische Universit\"{a}t Dresden\\
Dresden\\
Germany
\newline {\em and}
\newline
Instytut Matematyki\\
Uniwersytet Warszawski\\
Warszawa\\
Poland}
\email{daniel.max.hoffmann@gmail.com}
\urladdr{{https://sites.google.com/site/danielmaxhoffmann/}}

\author[P. KOWALSKI]{Piotr Kowalski$^{\diamondsuit}$}
\thanks{$^{\diamondsuit}$
 Supported by the Narodowe Centrum Nauki grant no. 2021/43/B/ST1/00405}
\address{$^{\diamondsuit}$Instytut Matematyczny\\
Uniwersytet Wroc{\l}awski\\
Wroc{\l}aw\\
Poland}
\email{pkowa@math.uni.wroc.pl} \urladdr{http://www.math.uni.wroc.pl/\textasciitilde pkowa/ }

\author[C.M. Tran]{Chieu-Minh Tran$^{\heartsuit}$}
\address{$^{\heartsuit}$Department of Mathematics \\ National University of Singapore \\ Singapore}
\email{trancm@nus.edu.sg}
\urladdr{https://blog.nus.edu.sg/tranchieuminhchieutran/}

\author[J. Ye]{Jinhe Ye$^{\spadesuit}$}
\address{$^{\spadesuit}$Mathematical Institute, University of Oxford}
\email{jinhe.ye@maths.ox.ac.uk}
\urladdr{https://sites.google.com/view/vincentye}

\date{\today}

 \DeclareMathOperator{\gal}{Gal}

 \DeclareMathOperator{\alg}{alg}

\DeclareMathOperator{\ddf}{DF}\DeclareMathOperator{\dcf}{DCF}\DeclareMathOperator{\scf}{SCF}

\DeclareMathOperator{\Nn}{\mathbb{N}}

\DeclareMathOperator{\Zz}{\mathbb{Z}}
\DeclareMathOperator{\Rr}{\mathbb{R}}
\DeclareMathOperator{\Qq}{\mathbb{Q}}
\DeclareMathOperator{\Cc}{\mathbb{C}}

\newcommand{\cupdot}{\mathbin{\mathaccent\cdot\cup}}

\newtheorem{theorem}{Theorem}[section]
\newtheorem*{theorem*}{Main Theorem}

\newtheorem{prop}[theorem]{Proposition}
\newtheorem{lemma}[theorem]{Lemma}
\newtheorem{cor}[theorem]{Corollary}
\newtheorem{fact}[theorem]{Fact}
\theoremstyle{definition}
\newtheorem{definition}[theorem]{Definition}
\newtheorem{example}[theorem]{Example}
\newtheorem{remark}[theorem]{Remark}

\newtheorem{notation}[theorem]{Notation}
\newtheorem{ass}[theorem]{Assumption}

\theoremstyle{remark}

\newtheorem*{cor*}{Corollary}

\theoremstyle{definition}
\theoremstyle{definition}

\theoremstyle{definition}

\theoremstyle{remark}

\makeatletter
\providecommand*{\cupdot}{%
  \mathbin{%
    \mathpalette\@cupdot{}%
  }%
}
\newcommand*{\@cupdot}[2]{%
  \ooalign{%
    $\m@th#1\cup$\cr
    \sbox0{$#1\cup$}%
    \dimen@=\ht0 %
    \sbox0{$\m@th#1\cdot$}%
    \advance\dimen@ by -\ht0 %
    \dimen@=.5\dimen@
    \hidewidth\raise\dimen@\box0\hidewidth
  }%
}

\providecommand*{\bigcupdot}{%
  \mathop{%
    \vphantom{\bigcup}%
    \mathpalette\@bigcupdot{}%
  }%
}
\newcommand*{\@bigcupdot}[2]{%
  \ooalign{%
    $\m@th#1\bigcup$\cr
    \sbox0{$#1\bigcup$}%
    \dimen@=\ht0 %
    \advance\dimen@ by -\dp0 %
    \sbox0{\scalebox{2}{$\m@th#1\cdot$}}%
    \advance\dimen@ by -\ht0 %
    \dimen@=.5\dimen@
    \hidewidth\raise\dimen@\box0\hidewidth
  }%
}
\makeatother

\def\Ind#1#2{#1\setbox0=\hbox{$#1x$}\kern\wd0\hbox to 0pt{\hss$#1\mid$\hss}
\lower.9\ht0\hbox to 0pt{\hss$#1\smile$\hss}\kern\wd0}

\def\notind#1#2{#1\setbox0=\hbox{$#1x$}\kern\wd0
\hbox to 0pt{\mathchardef\nn=12854\hss$#1\nn$\kern1.4\wd0\hss}
\hbox to 0pt{\hss$#1\mid$\hss}\lower.9\ht0 \hbox to 0pt{\hss$#1\smile$\hss}\kern\wd0}

\begin{document}

\newcommand{\ov}{\overline}
\newcommand{\FC}{\mathfrak{C}}

\newcommand{\twoc}[3]{ {#1} \choose {{#2}|{#3}}}
\newcommand{\thrc}[4]{ {#1} \choose {{#2}|{#3}|{#4}}}
\newcommand{\Kk}{{\mathds{K}}}

\newcommand{\dlog}{\mathrm{ld}}
\newcommand{\ga}{\mathbb{G}_{\rm{a}}}
\newcommand{\gm}{\mathbb{G}_{\rm{m}}}
\newcommand{\gaf}{\widehat{\mathbb{G}}_{\rm{a}}}
\newcommand{\gmf}{\widehat{\mathbb{G}}_{\rm{m}}}
\newcommand{\gdf}{\mathfrak{g}-\ddf}
\newcommand{\gdcf}{\mathfrak{g}-\dcf}
\newcommand{\fdf}{F-\ddf}
\newcommand{\fdcf}{F-\dcf}
\newcommand{\mw}{\scf_{\text{MW},e}}

\newcommand{\BC}{{\mathbb C}}

\newcommand{\CC}{{\mathcal C}}
\newcommand{\CG}{{\mathcal G}}
\newcommand{\CK}{{\mathcal K}}
\newcommand{\CL}{{\mathcal L}}
\newcommand{\CN}{{\mathcal N}}
\newcommand{\CS}{{\mathcal S}}
\newcommand{\CU}{{\mathcal U}}
\newcommand{\CF}{{\mathcal F}}
\newcommand{\CP}{{\mathcal P}}
\newcommand{\CI}{{\mathcal I}}
\newcommand{\SL}{{\mathrm{SL}}}

\begin{abstract}
    We show that if $G$ is a split semisimple algebraic group over a model complete field $K$,
    then the groups $G(K)$ and $G(K)'$ (the commutator group which is a ``Chevalley group'' as for example the group $\mathrm{PSL}_2(K)$) are model complete as well.
\end{abstract}

\maketitle

\section{Introduction}
The goal of this project is to show model completeness of some groups with geometric origin. Recall that
model completeness is a weaker and, hence, more flexible variant of quantifier elimination - we can reduce formulas to the form, where only a tuple of existential quantifiers remains. Many structures from algebra, geometry, and number theory are model complete but do not have quantifier elimination in their natural languages. Examples include fields\ $\mathbb
{R}$ of real numbers~\cite[Theorem 2.7.3]{HoMo}, the field $\mathbb{Q}_p$ of $p$-adic numbers~\cite{Mac76}, perfect PAC fields satisfying certain Galois-theoretic conditions~\cite{mcefree}, the exponential field~$(\mathbb{R}, \text{exp})$ of real numbers~\cite{Wilkie}, etc.

Another favorite source of examples for model theorists comes from the model companion construction, which serves as an analogue to the universal domain à la Weil for more delicate theories other than fields. Note that model companions are automatically model complete but usually do not have quantifier elimination. Examples include real closed fields (RCF)~\cite[Theorem 2.7.3]{HoMo}, $p$-adically closed fields (pCF)~\cite{Mac76}, algebraically closed valued fields (ACVF)~\cite{lou-notes}, algebraically closed fields with a  generic automorphism ACFA~\cite{acfa1}, or more exotic examples like $G$-TCF~\cite{OzlemnPiotr, HK3}, ACFO and interpolative fusions~\cite{minh, KTW}, ACFG~\cite{Christian},  CXFs and VXFs~\cite{cxf,vxf}.

An {\it algebraic group}, for us, is a smooth group scheme of finite type over a field $K$. An algebraic group $G$ is \emph{semisimple} if it is connected and any normal commutative subgroup of $G(K^{\alg})$ is finite.
Our main result (see Theorem \ref{sscomm}) is as follows (the notion of a \emph{split} algebraic group will be briefly discussed in Section \ref{secclass}):

\begin{theorem*}
Let $K$ be a model complete field and $G$ a semisimple and split algebraic group over $K$.
Then the groups $G(K)$ and $G(K)'$ (the commutator group) are model complete.
\end{theorem*}

In particular, our result implies both $\SL_2(\mathbb{R})$ and  $\SL_2(\mathbb{C})$ are model complete.
Initially, the project started without P.K., and we concentrated on these examples to demonstrate the feasibility of such results.
P.K. joined the project following D.M.H.'s talk at the Antalya Algebra Days conference in May 2023, contributing new ideas to establish model completeness for a broader class of groups.

There are several further motivations for the main theorem. For the field $\mathbb
{R}$ and the exponential field $(\mathbb
{R}, \text{exp})$, the completeness of the model can be used to show the geometric properties of these structures (more precisely, these structures are o-minimal, so definable sets and definable groups are very close to classical geometric objects like smooth manifolds and Lie groups). Our heuristic is that the implication also often goes the other way in the presence of groups, in particular, non-abelian groups of geometric significance such as $\SL_2(\mathbb{R})$ and  $\SL_2(\mathbb{C})$ should be model-complete. See also~\cite{Po2} for more on this philosophy. This view is also informed by a separated project on locally compact groups and Lie groups by D.M.H., where a related first-order structure is shown to have a good behavior provided the theory of the corresponding pure group is model complete.

Let us briefly shift the focus to the related problem of describing canonical topologies on natural examples in model theory (e.g., the Euclidean topology on the field of real numbers, the valuation topology on ACVF, etc). In the case of pure fields, the natural topology can be recovered by considering \'etale images, as outlined in~\cite{firstpaper}, with subsequent studies further expanding on this concept~\cite{field-top-1, field-top-2, ez}. This approach has limitations; for instance, in a model of ACVF, it incorrectly yields the Zariski topology,  where we do not have constructibility of definable sets. A viewpoint by C.M.T., also one of the authors of~\cite{firstpaper}, is that one would ultimately like to recover topology from group structures instead of fields. For $\SL_2(\mathbb{R})$ and $\SL_2(\mathbb{C})$, this should come from a suitable variant of the notion of \'etale images. Establishing the model completeness of $\SL_2(\mathbb{R})$ and  $\SL_2(\mathbb{C})$ is the first step in this direction.

Any interpretability functor takes elementary embeddings into elementary embeddings. However, such a functor \emph{need not} preserve the model completeness of structures (for example, $\Qq_p$ is a model complete field, but the group $\Qq_p^*$ is not model complete). In this paper, we consider a particular case when the interpretability functor comes either from a semisimple algebraic group or the corresponding ``Chevalley group''. It is worth noticing that (to our knowledge) not many examples of infinite model complete groups are known besides the commutative ones.
Nick Ramsey have pointed out to us that the only known examples of non-commutative infinite model complete groups are extraspecial groups (see e.g. \cite[Theorem 3.2]{baudneo}) which are still 2-nilpotent, so they are ``very close'' to commutative groups.
Therefore, our paper gives in particular first examples of model complete infinite simple groups.

In the previous version of this paper, interpreting a field within an algebraic group was used in an essential way. We avoid it in this version using a model-theoretic trick suggested to us by Will Johnson (see Theorem \ref{will}). Still, such interpretations continue to be an interesting research direction in model theory starting from the work of Simon Thomas in the case of ``Chevalley group'' (see \cite{thomasthesis}). Later, Ali Nesin achieved this for $\text{SO}_3(\Rr)$ in~\cite{NESIN1989}. A similar result is obtained in \cite{NesinPillay91}, in Theorem 1.1 from \cite{PePiSe}, and more recently in \cite{segtent}.

To highlight our result,
in Theorem \ref{finalprod}, we show that homomorphisms between simply connected semisimple algebraic groups decompose into a field homomorphism part and an algebraic isogeny part. Proposition \ref{mainmulti} provides the elementarity of group maps induced by several field maps. Finally, Theorem \ref{subgroup} gives a result about extensions of automorphisms of finite commutative groups. Having all of that on board, we can prove the Main Theorem (Theorem \ref{sscomm}) about the model completeness of the groups $G(K)$ and $G(K)'$, where $G$ is a split semisimple algebraic group over a model complete field $K$. In this proof, we treat the (much easier) ``Chevalley case'' $G(K)'$ first and then use it (as well as some other methods) to show the model completeness of $G(K)$.

The main text is divided into two sections. In Section 2, we recall definitions and classical facts from model theory and the theory of algebraic groups, as well as prove the necessary intermediate results outlined above. In Section 3 we obtain the main result about model completeness of groups coming from semisimple split algebraic groups.

The authors would like to thank  Zo\'e Chatzidakis, Philip Dittmann, Lou van den Dries, Martin Hils, Aleksander Ivanov, Franziska Jahnke,  Will Johnson, Nick Ramsey, Tomasz Rzepecki, Katrin Tent, and Martin Ziegler  for various discussions on this paper.
As each of the authors of this paper works in a different place, the work was conducted during meetings in many places (Beijing, Luminy, M\"{u}nster, Singapore, Sirince, South Bend, Wrocław).
We would like to thank the members of the
model theory group in Wroc{\l}aw for their constructive remarks during the talks of P.K. at the model theory
seminars at Wroc{\l}aw University.

\section{Homomorphism between rational points of algebraic groups}
In this section we collect the definitions and results needed for the statement and the proof of our main results. In Subsection \ref{secclass}, some classical results about rational points of semisimple algebraic groups are collected culminating at the Borel-Tits Theorem. In Subsection \ref{secconbt}, we show some consequences of the Borel-Tits Theorem regarding homomorphisms between rational points of semisimple algebraic groups. In Subsection \ref{secmt}, we provide a few notions and facts from model theory. Subsection \ref{secsf} is crucial: we consider there a ``multi-field situation'' and show that certain group homomorphisms induced by several field homomorphisms are elementary. Subsection \ref{secaut} seems not to fit much to the topic of this paper, however we use the main result of it (Theorem \ref{subgroup}) for showing the model completeness in the semisimple case.

\subsection{Classical results on algebraic groups}\label{secclass}
We mostly follow Milne's book \cite{milnealggps}. It is often not very important to provide a given definition fully, when we use particular consequences of this definition only.

All the algebraic groups we consider are over fields and they are smooth and connected. For an algebraic group $G$ over $K$ and a field extension $K\subseteq L$, by $G_L$ we mean the change of scalars from $K$ to $L$. By a \emph{simple algebraic group} over a field $K$, we mean an infinite algebraic group $G$ over $K$ such that any proper normal subgroup of $G(K^{\alg})$ is finite.

By the \emph{radical} of an algebraic group $G$ over $K$, denoted $\mathrm{Rad}(G)$, we mean its maximal connected normal algebraic solvable subgroup over $K$. By a \emph{semisimple algebraic group} $G$ over a field $K$, we mean an infinite algebraic group over $K$ such that $\mathrm{Rad}(G_{K^{\alg}})$ is trivial, equivalently, any normal commutative subgroup of $G(K^{\alg})$ is finite.

An \emph{isogeny} is an algebraic group epimorphism with a finite kernel. For the definition of a \emph{multiplicative} isogeny, we refer the reader to \cite[Definition 18.1]{milnealggps}. A semisimple algebraic group $G$ is \emph{simply connected} if every
multiplicative isogeny $G_0\to G$ of connected algebraic groups is an isomorphism (see \cite[Definition 18.5]{milnealggps}). A semisimple algebraic group $G$ is \emph{adjoint} if $Z(G)$ is trivial.

For the notion of a \emph{split} algebraic group, we refer the reader to \cite[Definition 17.101]{milnealggps}. We will need it only for Theorems \ref{need} and \ref{stein}, which we will be using as black boxes anyway.
\begin{remark}\label{term}
For other possible choices of the terminology above, the reader is advised to consult \cite[Definition 19.8]{milnealggps}. General definitions are given e.g. here \cite[Def XIX, 2.7]{sga33}.
\end{remark}

We state below several classical results.
\begin{theorem}[Chapter 24a in \cite{milnealggps}]\label{prod}
If $G$ is a semisimple and simply connected algebraic group over a field $K$, then there is a decomposition
      $$G\cong G_1\times \ldots G_l$$
where $G_1,\ldots,G_l$ are simple and simply connected algebraic groups over $K$.
\end{theorem}
\noindent
This result might be found as well in \cite{margulis}, see Prop. 1.4.10 there. It is worth mentioning that this decomposition is unique up to isomorphism. We will use the following.
\begin{notation}\label{scnot}
Let $G$ be a semisimple algebraic group over a field $K$.  Then the algebraic subgroup $Z(G)$ is finite and the semisimple algebraic group $G/Z(G)$ is adjoint. We denote:
$$G_{\mathrm{ad}}:=G/Z(G)$$
 and the corresponding quotient isogeny by
 $$\pi_{\mathrm{ad}}:G\to G_{\mathrm{ad}}.$$
 There is also a ``universal cover'' isogeny
 $$\pi_{\mathrm{sc}}:G_{\mathrm{sc}}\to G$$
 over $K$, where $G_{\mathrm{sc}}$ is a simply connected semisimple algebraic group over $K$ (see \cite[Chapter 18d]{milnealggps}).
\end{notation}
We will need the following results about rational points of semisimple groups.
\begin{theorem}\label{need}
Let $G$ be a semisimple split algebraic group over an infinite field $K$. Then we have the following.
\begin{enumerate}
\item Any proper normal subgroup of $G(K)'$ is contained in $Z(G(K))$.


\item If $G$ is simply connected, then $G(K)$ is perfect.

\item The group $G(K)'$ is perfect.


\item The group $G(K^{\alg})$ is perfect.

\item Suppose that $G$ is simple. Then $\ker(\pi_{\mathrm{sc}})\cong \mu_2\times \mu_2$ or there is $n>0$ such that $\ker(\pi_{\mathrm{sc}})\cong \mu_n$, where $\mu_n$ is the group scheme of $n$-th roots of unity.
    
\end{enumerate}
\end{theorem}
\begin{proof}
By \cite[Corollaire 6.5]{Botits} and \cite[Remarque 6.6]{Botits} (since $G$ is split), we have $G(K)'=G^+$ (or $G^+_K$ in the terminology from \cite[Section 7.2]{plat93}). Therefore, we can use \cite[Theorem 7.1]{plat93} and we obtain Item $(1)$.
\\
Item $(2)$ follows again from \cite[Corollaire 6.5]{Botits} and \cite[Remarque 6.6]{Botits} (since $G$ is split).
\\
Item $(3)$ follows as above using also \cite[Corollaire 6.4]{Botits}.
\\
Item $(4)$ holds even in a stronger form: every element of $G(K^{\alg})$ is a commutator, see \cite[21.58]{milnealggps}.
\\
For Item $(5)$, see Table 9.2 on page 71 of \cite{linlie}.
\end{proof}

We immediately obtain.
\begin{cor}\label{homchev}
Let $G,K$ be as above, $W$ be an arbitrary group, and $\alpha:G(K)'\to W$ be a homomorphism. Then, the following are equivalent.
\begin{enumerate}
  \item $\alpha$ is non-trivial.
  \item $\ker(\alpha)$ is finite and central.
  \item $\mathrm{im}(\alpha)$ is infinite.
\end{enumerate}
\end{cor}
The result below was proved by Rosenlicht in the case of a perfect field $M$ and by Grothendieck in the arbitrary case. We will need it for model complete fields which are necessarily perfect.
\begin{theorem}[\cite{ros57} and \cite{sga33}]\label{zd}
If $G$ is a reductive algebraic group defined over an infinite field $K$ (in particular: it applies to semisimple algebraic groups), then $G(K)$ is Zariski dense in $G(K^{\alg})$.
\end{theorem}
For a homomorphism of algebraic groups $f:H\to G$ over $K$, we denote by $f_K:H(K)\to G(K)$ the corresponding group homomorphism between the rational points. If $\varphi:K\to K'$ is a field homomorphism, then $\varphi_G:G(K)\to {}^{\varphi}G(K')$ denotes the corresponding group homomorphism between the rational points, where ${}^{\varphi}G$ denotes the corresponding (change of basis) algebraic group over $K'$ (we have $G(K')={}^{\varphi}G(K')$). We will often use the following crucial result by Borel and Tits \cite[Theorem 1.3]{Steinberg} (the original reference is \cite[(A)]{Botits}), which we state in a slightly simplified form.
\begin{theorem}[Theorem 1.3 in \cite{Steinberg}]\label{stein}
Let $H,G$ be simple split algebraic groups defined over infinite fields $L,M$ respectively. Assume that $H$ is simply connected or $G$ is adjoint. Let $W$ be a subgroup of $H(L)$ containing $H(L)'$ and $\alpha:W\to G(M)$ be a group homomorphism such that $\alpha(W)$ is Zariski dense. Then there are unique:
\begin{itemize}
  \item a field homomorphism $\varphi :L\to M$;
  \item a ``special''  (see Remark \ref{zdense}(3)) isogeny $\beta:{}^{\varphi}H\to G$
\end{itemize}
such that we have:
$$\alpha=\beta_M\circ \varphi_H|_{W}.$$
\end{theorem}
\begin{remark}\label{zdense}
We comment here on several issues related to the Borel-Tits Theorem above.
\begin{enumerate}
  \item By Theorem \ref{zd} the assumption that $\alpha(H(L))$ is Zariski dense (in the statement of Theorem \ref{stein}) is unambiguous since being Zariski dense in $G(M)$ is the same as being Zariski dense in $G(M^{\alg})$.

  \item By \cite[Remark 1.4(b)]{Steinberg}, Theorem \ref{stein} is indeed a special case of \cite[Theorem 1.3]{Steinberg}.

  \item The isogeny $\beta$ above is ``special'' if $\mathrm{d}\beta\neq 0$ (see \cite[Theorem 1.3]{Steinberg}). We will not comment more on this condition, since it is needed for the uniqueness claim only (in the statement of Theorem \ref{stein}).
\end{enumerate}
\end{remark}

\subsection{Consequences of Borel-Tits Theorem}\label{secconbt}
In this subsection, we collect some consequences of Theorem \ref{stein}.

\begin{prop}\label{mainlemma}
Assume that $H$ is a split simple algebraic group over an infinite field $L$, $G$ is a connected linear algebraic group over a field $M$, and $\alpha:H(L)'\to G(M)$ is a non-trivial homomorphism. Then we have the following:
\begin{enumerate}
  \item  $\dim(H)\leqslant \dim(G)$;
  \item if $\dim(H)=\dim(G)$, then the image of $\alpha$ is Zariski dense.
\end{enumerate}
\end{prop}
\begin{proof}
Without loss, we can assume $M=M^{\alg}$. By composing with the ``universal cover'' isogeny $\pi_{\mathrm{sc}}:H_{\mathrm{sc}}\to H$, we can also assume that $H$ is simply connected. By Corollary \ref{homchev}, $\alpha $ has infinite image. Let $G_0$ denote the algebraic subgroup of $G$ being the Zariski closure of $\alpha(H(L))$ in $G(M^{\alg})$ (we identify here $G_0$ with $G_0(M^{\alg})$).

Since $H$ is semisimple, we obtain that $\mathrm{Rad}(G_0)(M)\cap \mathrm{im}(\alpha)$ is finite. Let us define:
$$G_1:=(G_0)_{\mathrm{ss}}:=G_0/\mathrm{Rad}(G_0).$$
By \cite[Proposition 19.2(b)]{milnealggps}, $G_1$ is semisimple. Then, the following composition map:
$$H(L)\longrightarrow G_0\left(M^{\alg}\right)\longrightarrow G_1\left(M^{\alg}\right)$$
has infinite image that is dense.
Since $G_1$ is semisimple, by Theorem \ref{prod}(1), the existence of ``universal cover'' and the same argument as above, we can assume that $G_1$ is simple. We apply now Theorem \ref{stein} for the induced homomorphism $H(L)\to G_1(M^{\alg})$, hence we obtain a field homomorphism $\varphi: L\to M^{\alg}$ and an isogeny ${}^{\varphi}H\to G_1$. Therefore we get:
$$\dim(H)=\dim\left({}^{\varphi}H\right)=\dim(G_1)\leqslant \dim(G)$$
giving Item $(1)$.

In particular, $\dim(H)=\dim(G)$ implies that $G_1=G_0=G$, hence the image of $\alpha$ is Zariski dense which is Item $(2)$.
\end{proof}
\begin{remark}
\begin{enumerate}
\item The algebraic group $G_0$ appearing in the statement of Proposition \ref{mainlemma} need not be simple, for example let us consider: $$H=\mathrm{SL}_2,\ G=\mathrm{SL}_2\times \mathrm{SL}_2,\ L=M=\Cc,\ \alpha(g)=(g,\bar{g}),$$
    where $g\mapsto \bar{g}$ is the complex conjugation.


  \item By the Chevalley's structure theorem, we could have dropped the linearity assumption on $G$.
\end{enumerate}
\end{remark}
We will need a result about homomorphisms between rational points of products of simple group schemes.
\begin{theorem}\label{finalprod}
Let $S_i,T_i$ be simply connected simple and split group schemes defined over infinite fields $L,M$ respectively such that $\dim(S_i)=\dim(T_i)$ ($i=1,\ldots,l$). Let $N$ be an algebraic subgroup of $Z(T_1\times \ldots \times T_l)$, $H:=(T_1\times \ldots \times T_l)/N$ and
 $$\Psi:S_1(L)\times \ldots \times S_l(L)\longrightarrow H(M)$$
be a group homomorphism with finite kernel. Then, there are $\varphi_i:L\to M$ and an isogeny
$$\beta:{}^{\varphi_1}S_1\times \ldots \times {}^{\varphi_l}S_l \to H$$
such that
$$\Psi=\beta_M\circ\left((\varphi_1)_{S_1}\times \ldots \times (\varphi_l)_{S_l}\right).$$
\end{theorem}
\begin{proof}
For simplicity of the presentation, we assume that $l=2$. We consider the quotient isogeny (see Notation \ref{scnot}):
$$\pi_{\mathrm{ad}}:H\to H_{\mathrm{ad}}=(T_1)_{\mathrm{ad}}\times (T_2)_{\mathrm{ad}},$$
and the corresponding homomorphism with finite kernel:
$$\widetilde{\Psi}:S_1(L)\times S_2(L)\to (T_1)_{\mathrm{ad}}(M)\times (T_2)_{\mathrm{ad}}(M),\ \ \ \widetilde{\Psi}=\pi\circ \Psi.$$
{\bf Claim 1}
\\
There is $\tau\in \mathrm{Sym}(\{1,2\})$ and homomorphisms ($i=1,2$)
 $$ \Psi_i:S_i(L)\longrightarrow (T_{\tau(i)})_{\mathrm{ad}}(M)$$
 such that $\widetilde{\Psi}=\tau\circ(\Psi_1\times \Psi_2)$.
\begin{proof}[Proof of Claim 1]
We have four homomorphisms $\Psi^i_j:S_i(L)\to (T_j)_{\mathrm{ad}}(M)$ for $i,j\in \{1,2\}$. For example, $\Psi^1_1$ is defined as the following composition:
$$S_1(L)\hookrightarrow S_1(L)\times S_2(L)\xrightarrow{\widetilde{\Psi}} (T_1)_{\mathrm{ad}}\times (T_2)_{\mathrm{ad}} \twoheadrightarrow (T_1)_{\mathrm{ad}}.$$
We need to show that $\Psi^1_2,\Psi^2_1$ are trivial or $\Psi^1_1,\Psi^2_2$ are trivial.
\\
\\
{\bf Subclaim}
\\
If $\Psi^1_1$ is non-trivial, then $\Psi^2_1$ is trivial (similarly for $\Psi^2_2$ and $\Psi^1_2$).
\begin{proof}[Proof of Subclaim]
By Proposition \ref{mainlemma} (using $\dim(S_i)=\dim(T_i)$), the image of $\Psi^1_1$ is Zariski dense. Therefore, we have:
$$C_{(T_1)_{\mathrm{ad}}(M)}\left(\Psi^1_1\left(S_1(L)\right)\right)=Z\left((T_1)_{\mathrm{ad}}(M)\right).$$
Since
$$\Psi^2_1\left(S_2(L)\right)\subseteq C_{(T_1)_{\mathrm{ad}}(M)}\left(\Psi^1_1\left(S_1(L)\right)\right),$$
we obtain that the image of $\Psi^2_1$ is a commutative group. Since $S_2(L)$ is perfect (by Theorem \ref{need}(2)), $\Psi^2_1$ is trivial.
\end{proof}
We consider two cases.
\\
\\
{\bf Case 1} $\Psi^1_1$ is non-trivial.
\\
By Subclaim, $\Psi^2_1$ is trivial. Since $\Psi$ is one-to-one, we obtain that $\Psi^2_2$ is non-trivial. Using Subclaim again, we get that $\Psi^1_2$ is trivial.
\\
\\
{\bf Case 2} $\Psi^1_1$ is trivial.
\\
Since $\Psi$ is one-to-one, we obtain that $\Psi^1_2$ is non-trivial. By Subclaim, we get that $\Psi^2_2$ is trivial.
\end{proof}
By Proposition \ref{mainlemma}(1) and Claim, we get that
$$\dim(S_i)=\dim\left(T_{\tau(i)}\right)$$
for $i=1,2$. Hence we can assume for simplicity of the presentation that $\tau=\mathrm{id}$ (in the general case, we would need to write ``$\tau(i)$'' in several places and compose at the end with the algebraic isomorphism given by the permutation of coordinates corresponding to $\tau$) and denote
$$\widetilde{T_i}:=T_i/(N\cap T_i)\leqslant H.$$
Then we have $[\widetilde{T_1},\widetilde{T_2}]=1$ and the multiplication map $m:\widetilde{T_1}\times \widetilde{T_2}\to H$ is an isogeny.
\\
\\
{\bf Claim 2}
\\
$\Psi(S_i(L))\subseteq \widetilde{T_i}(M)$.
\begin{proof}[Proof of Claim 2]
Using Claim 1, we obtain:
$$\Psi(S_i(L))\subseteq \widetilde{T_i}(M)Z(H(M)).$$
Let us consider the following composition map $\alpha_i$:
$$S_i(L)\xrightarrow{\Psi} \widetilde{T_i}(M)Z(H(M))\to \widetilde{T_i}(M)Z(H(M))/\widetilde{T_i}(M) \cong Z(H(M))/\left(\widetilde{T_i}(M)\cap Z(H(M))\right).$$
Since $S_i(L)=S_i(L)'$ and $Z(H(M))/\left(\widetilde{T_i}(M)\cap Z(H(M))\right)$ is commutative, we obtain that the map $\alpha_i$ is trivial and $\Psi(S_i(L))\subseteq \widetilde{T_i}(M)$.
\end{proof}
Let us denote:
$$\Psi_i:=\Psi|_{S_i(L)}:S_i(L)\to \widetilde{T_i}(M).$$
Using our dimension assumptions and Proposition \ref{mainlemma}(2), we apply Theorem \ref{stein} to $\Psi_i$ to obtain $\varphi_i :L\to M$ and isogenies $\beta_i:{}^{\varphi}S_i\to \widetilde{T_i}$
such that:
$$\Psi_i=(\beta_i)_M\circ (\varphi_i)_{S_i}.$$
Then it is enough to take $\beta:=m\circ (\beta_1\times \beta_2)$.
\end{proof}
\begin{cor}\label{ssdense}
If $G$ (over an infinite field $L$) and $H$ (over an infinite field $M$) are semisimple such that
$$G_{\mathrm{sc}}=S_1\times \ldots \times S_n,\ \ H_{\mathrm{sc}}=T_1\times \ldots \times T_n,$$
where $S_i,T_i$ are simple algebraic and such that $\dim(T_i)=\dim(S_i)$. If $f:G(L)'\to H(M)$ is a homomorphism with finite kernel, then $\mathrm{im}(f)$ is Zariski dense. In particular, we have
$$f(Z(G(L)'))\subseteq Z(H(M)).$$
\end{cor}
\begin{proof}
Let $\widetilde{f}:=f\circ \pi$, where $\pi:G_{\mathrm{sc}}(L)\to G(L)$ is the ``universal cover'' map. Then $\widetilde{f}$ has still finite kernel, so we can assume that $G$ is simply connected. We are now in the situation of Theorem \ref{finalprod}, so we get that
$$f=\beta_M\circ\left((\varphi_1)_{S_1}\times \ldots \times (\varphi_l)_{S_l}\right),$$
where $\beta_M$ is an isogeny and each $(\varphi_i)_{S_i}$ has Zariski dense image by Theorem \ref{zd}. Since the product of Zariski dense subsets is Zariski dense (also with respect to the Zariski topology on the product, which is richer than the product topology!) and isogenies are dominant, we get that image of $f$ is Zariski dense.

The ``in particular'' claim follows.
\end{proof}

We need one more consequence of the Borel-Tits Theorem.
\begin{lemma}\label{coincide}
Suppose that $G$ is a simple adjoint algebraic group over an infinite field $K$. Let $L,M$ be field extensions of $M$, $G(L)'\leqslant W\leqslant G(L)$, and $f_1,f_2:W\to G(M)$ be group homomorphisms with finite kernels and agreeing on $G(L)'$. Then $f_1=f_2$.
\end{lemma}
\begin{proof}
By Proposition \ref{mainlemma}, $f_1(G(L)')$ and $f_2(G(L)')$ are Zariski dense in $G(M)$. By Theorem \ref{stein}, there are field maps $\varphi,\varphi_1,\varphi_2:L\to M$ and ``special'' isogenies (see Remark \ref{zdense}(3))
$$\beta:{}^{\varphi}G\to G,\ \ \beta_1:{}^{\varphi_1}G\to G,\ \ \beta_2:{}^{\varphi_2}G\to G$$
such that
$$f_1|_{G(L)'}=\beta_M\circ \varphi_G|_{G(L)'}=f_2|_{G(L)'},\ \ \ f_1=(\beta_1)_M\circ (\varphi_1)_G|_{W},\ \ \
f_2=(\beta_2)_M\circ (\varphi_2)_G|_{W}.$$
By the uniqueness claim in Theorem \ref{stein} applied to the restrictions to $G(L)'$, we get:
$$\varphi_1=\varphi=\varphi_2,\ \ \ \beta_1=\beta=\beta_2$$
and  the result follows.
\end{proof}
\subsection{Miscellaneous model theory}\label{secmt}
We provide here some notions and facts from model theory.
\begin{definition}\label{mcparam}
Let $\mathcal{L}$ be a language and $M$ be an $\mathcal{L}$-structure. We say that $M$ is \emph{model complete} if $\mathrm{Th}(M)$ is
    model complete.
\end{definition}
\begin{remark}\label{mcspecial}
To test model completeness of a theory $T$, it is enough to consider monomorphisms between ``sufficiently nice'' (e.g. special, resplendent, $\kappa$-saturated for a convenient cardinal $\kappa$) models of $T$.
\end{remark}
The following approach was suggested to us by Will Johnson. Let $T_i$ is an $\mathcal{L}_i$-theory for $i=1,2$ and $\Gamma:\mathrm{Models}(T_1)\to \mathrm{Models}(T_2)$ be an interpretability functor (categories of models with elementary embeddings as morphisms) as in \cite[Theorem 5.3.4]{HoMo}.
\begin{example}
If $G$ is an algebraic group over $K$, then we can take:
\begin{itemize}
  \item $\mathcal{L}_1:=(+,\cdot,-,0,1,\lambda_a)_{a\in K}$ as the language of $K$-algebras;
  \item $T_1$ as the $\mathcal{L}_1$-theory of field extensions of $K$;
  \item $\mathcal{L}_2$ as the language of groups;
  \item $T_2$ as the theory of groups.
\end{itemize}
Then
$$\Gamma:\mathrm{Models}(T_1)\to \mathrm{Models}(T_2),\ \ \ \ \Gamma(K)=G(K)$$
is an interpretability functor.
\end{example}
We will need the following consequence of Robinson's joint consistency theorem (see \cite[Corollary 9.5.8]{HoMo}).
\begin{theorem}\label{will}
Suppose that $\Gamma:\mathrm{Models}(T_1)\to \mathrm{Models}(T_2)$ is an interpretability functor
and  $M_2\models T_2$ is special (as in \cite[Section 10.4]{HoMo}).
If there is $M_1\models T$ such that $M_2\equiv \Gamma(M_1)$, then there is $M_1'\models T$ such that $M_2\cong \Gamma(M_1')$.
\end{theorem}
\begin{proof}
Let $\mathcal{L}$ be the two-sorted language $\mathcal{L}_1\cup \mathcal{L}_2\cup \{f\}$, where $f$ is the a function symbol from $\mathcal{L}_1$-structures to $\mathcal{L}_2$-structures. Let $T$ be the $\mathcal{L}$-theory expressing that the $\mathcal{L}_1$-sort is a model of $T_1$, the $\mathcal{L}_2$-sort is a model of $T_2$ and $f$ is an $\mathcal{L}_2$-isomorphism between $\Gamma(M_1)$ and the $\mathcal{L}_1$-sort.

By \cite[Prop. 1.2(2)]{CasExp}, $M_2$ can be expanded to an $\mathcal{L}$-structure which is a model of $T$, so we obtain $M_1'$ as in the statement.
\end{proof}
The following result will not be directly used but we include it as an illustration.
\begin{cor}\label{ziegler}
Let $M_1,M_2,M$ be $\mathcal{L}$-structures, where $\mathcal{L}$ is an arbitrary language. Suppose that $M\equiv M_1\times M_2$ and $M$ is special. Then there are $\mathcal{L}$-structures $N_1,N_2$ such that
$$M_1\equiv N_1,\ \ \ M_2\equiv N_2,\ \ \ M\cong N_1\times N_2.$$
\end{cor}
\begin{proof}
We use Theorem \ref{will} for $\mathcal{L}_1:=$``bi-sorted $\mathcal{L}$'',  $\mathcal{L}_2:=\mathcal{L}$, and $\Gamma$ being the product.
\end{proof}
\begin{remark}
We comment here on the product case and include one warning
\begin{enumerate}
  \item The result above is not true for arbitrary structures. As an easy example, we take:
  $$(\Qq,+)\equiv (\Qq,+)\times (\Qq,+),$$
  but $(\Qq,+)$ is not isomorphic to any non-trivial product.

  \item This result with two proofs (different than the one above) was pointed out to us by Martin Ziegler.

  \item Martin Hils pointed out to us an interesting example which we would like to include here. The model completeness is \emph{not} preserved by taking products even in the case of groups: the groups $C_{p^{\infty}}$ (the Pr\"{u}fer group) and $\Zz/2\Zz$ are model complete, but their product $C_{p^{\infty}}\times \Zz/2\Zz$ is not.


\end{enumerate}
\end{remark}
The following obvious observation will be used at the very end of this paper.
\begin{lemma}\label{obv}
Suppose that $f,g:N\to N'$ are monomorphisms of $\mathcal{L}$-structures such that $f(H)=g(H)$ and $g$ is elementary (resp. 1-elementary). Then $f$ is elementary (resp. 1-elementary).
\end{lemma}
\begin{proof}
It is clear, since a monomorphism is elementary if and only if its image is an elementary  (resp. 1-elementary) substructure.
\end{proof}

We finish this subsection with a special case of the Feferman-Vaught Theorem.
\begin{fact}[Cor 9.6.5(b) in \cite{HoMo}]\label{rem: Hodges2}
If $M_1\preccurlyeq N_1,M_2\preccurlyeq N_2$ then $M_1\times M_2\preccurlyeq N_1\times N_2$.
\end{fact}

\begin{remark}\label{ziegdisj}
Martin Ziegler pointed out to us that to show Fact \ref{rem: Hodges2} one only needs to consider the ``disjoint union'' (two-sorted) theory $\mathrm{Th}(M_1)\cupdot  \mathrm{Th}(M_2)$ rather than the (much stronger) Feferman-Vaught Theorem.
\end{remark}

\subsection{Several fields}\label{secsf}
We need some results about elementarity of homomorphisms between rational points of algebraic groups induced by \emph{several} elementary field maps.
\begin{ass}\label{assseveral}
All the fields considered in this subsection are perfect. Moreover, we assume the following.
\begin{enumerate}
  \item Let
  $$\varphi_1,\varphi_2:L\to M$$ be field homomorphisms extending to field homomorphisms
  $$\overline{\varphi}_1,\overline{\varphi}_2:L^{\alg}\to M^{\alg}$$
  respectively.

  \item Let $\widetilde{G}$ be a semisimple simply connected algebraic group over $L$ such that
 $$\widetilde{G}=S_1\times S_2,$$ where $S_1,S_2$ are simple algebraic groups over $L$.
 \item We denote:
 $$\widetilde{H}:={}^{\varphi_1}S_1\times {}^{\varphi_2}S_2,\ \ \ \ \ \ \
 \overline{\varphi}_{\widetilde{G}}:=\left(\overline{\varphi}_1\right)_{S_1}\times \left(\overline{\varphi}_2\right)_{S_2}:\widetilde{G}\left(L^{\alg}\right)\to \widetilde{H}\left(M^{\alg}\right).$$

 \item We have semisimple algebraic groups $G,H$ and the following isogenies (over $L$ and $M$ resp.):
 $$\alpha:\widetilde{G}\to G,\ \ \ \ \ \beta:\widetilde{H}\to H$$
 such that
 $$\overline{\varphi}_{\widetilde{G}}\left(N_{\alpha}\left(L^{\alg}\right)\right)=N_{\beta}\left(M^{\alg}\right),$$
 where $N_{\alpha}=\ker(\alpha)$ and $N_{\beta}=\ker(\beta)$.
\end{enumerate}
\end{ass}
From Assumption \ref{assseveral}(4), we get the induced map:
$$ \overline{\varphi}_{G}:G\left(L^{\alg}\right)\to H\left(M^{\alg}\right).$$
\begin{prop}\label{inducedg}
We have
$$ \overline{\varphi}_{G}(G\left(L\right))\subseteq H\left(M\right).$$
\end{prop}
\begin{proof}
We define:
$$\widetilde{G}_{\alpha}\left(L\right):=\left(\alpha_{L^{\alg}}\right)^{-1}(G(L))\leqslant \widetilde{G}\left(L^{\alg}\right).$$
\textbf{Claim 1
}
\\
We have:
$$\widetilde{G}_{\alpha}\left(L\right)=\left\{g\in \widetilde{G}\left(L^{\alg}\right)\ |\ \mathrm{Gal}(L)g\subseteq gN_{\alpha}\left(L^{\alg}\right)\right\},$$
where $\mathrm{Gal}(L):=\mathrm{Gal}\left(L^{\alg}/L\right)$ is the absolute Galois group of $L$.
\begin{proof}[Proof of Claim 1]
For any $g\in \widetilde{G}\left(L^{\alg}\right)$ and any $\sigma\in \mathrm{Gal}(L)$ we have that $\sigma_{\widetilde{G}}(g)\in gN_{\alpha}(L^{\alg})$ if and only if
$$\alpha_{L^{\alg}}(g)=\alpha_{L^{\alg}}(\sigma_{\widetilde{G}}(g))=\sigma_G\left(\alpha_{L^{\alg}}(g)\right).$$
Therefore, $\mathrm{Gal}(L)g\subseteq gN_{\alpha}(L^{\alg})$ if and only if
$$\alpha_{L^{\alg}}(g)\in G\left(L^{\alg}\right)^{\mathrm{Gal}(L)}=G(L),$$
which concludes the proof of Claim 1.
\end{proof}
We have also the following corresponding group
$$\widetilde{H}_{\beta}\left(M\right)\leqslant \widetilde{H}\left(M^{\alg}\right).$$
and (by the obvious analogue of Claim 1):
$$\widetilde{H}_{\beta}\left(M\right)=\left\{g\in \widetilde{H}\left(M^{\alg}\right)\ |\ \mathrm{Gal}(M)g\subseteq gN_{\beta}\left(L^{\alg}\right)\right\}.$$
\textbf{Claim 2
}
\\
We have the following:
$$\overline{\varphi}_{\widetilde{G}}\left(\widetilde{G}_{\alpha}\left(L\right)\right)\subseteq
\widetilde{H}_{\beta}\left(M\right).$$
\begin{proof}[Proof of Claim 2]
Let us take $\sigma\in \gal(M)$ and $i=1,2$. Since
$$\overline{\varphi}_i(L)=\varphi_i(L)\subseteq M,$$
$\sigma$ preserves $\overline{\varphi}_i(L)$ pointwise and then $\sigma$ preserves $\overline{\varphi}_i(L^{\alg})=\overline{\varphi}_i(L)^{\alg}$ setwise. Let $\tau_i$ denote the following composition map:
\begin{equation*}
 \xymatrix{L^{\alg} \ar[rr]^{\overline{\varphi}_i}  & &  \overline{\varphi}_i(L^{\alg})   \ar[rr]^{\sigma|_{\overline{\varphi}_i(L^{\alg})}}  & & \overline{\varphi}_i(L^{\alg})   \ar[rr]^{\overline{\varphi}_i^{-1}} & &   L^{\alg} .}
\end{equation*}
Then $\tau_i\in \gal(L)$ and we have:
$$\sigma\circ \overline{\varphi}_i=\overline{\varphi}_i\circ \tau_i.$$
Let us take now $g=(g_1,g_2)\in \widetilde{G}_{\alpha}(L)$ and compute:
\begin{IEEEeqnarray*}{rCl}
\sigma_{\widetilde{H}}\left(\overline{\varphi}_{\widetilde{G}}(g)\right) & = & \sigma_{\widetilde{H}}\left(\left(\overline{\varphi_1}\right)_{S_1}(g_1),\left(\overline{\varphi_2}\right)_{S_2}(g_2)\right) \\
 & = & \left(\left(\sigma\circ \overline{\varphi}_{1}\right)_{S_1}(g_1),\left(\sigma\circ \overline{\varphi}_{2}\right)_{S_2}(g_2)\right)\\
 & = & \left(\left(\overline{\varphi}_{1}\circ \tau_1\right)_{S_1}(g_1),\left(\overline{\varphi}_{2}\circ \tau_2\right)_{S_2}(g_2)\right)\\
 & = & \overline{\varphi}_{\widetilde{G}}\left(\left(\tau_1\right)_{S_1}(g_1),\left(\tau_2\right)_{S_2}(g_2)\right)\\
 &\in & \overline{\varphi}_{\mathrm{sc}}\left(gN_{\alpha}(L^{\alg})\right)\\
 &=& \overline{\varphi}_{\mathrm{sc}}(g)N_{\beta}(M^{\alg}),
\end{IEEEeqnarray*}
by Claim 1 and the last bullet point in Assumption \ref{assseveral} (and a rather trivial observation that the invariants of the diagonal action coincide with the invariants of the coordinate action). Using Claim 1 again, we get that $\overline{\varphi}_{\widetilde{G}}(g)\in \widetilde{H}_{\beta}\left(M\right)$, which finishes the proof of Claim 2.
\end{proof}
The result follows immediately using Claim 2.
\end{proof}
Let us denote the induced map from Proposition \ref{inducedg} by:
$$\varphi_{G}:G\left(L\right)\to H\left(M\right).$$
By Claim 1 from the proof of Proposition \ref{inducedg}, this map does not depend on the extensions $ \overline{\varphi}_i\supseteq \varphi_i$ (satisfying Assumption \ref{assseveral}).

In the next result we check the 1-elementarity only, since we aim to apply Robinson's test later in the proof of Theorem \ref{sscomm} in the case of $G(K)$. In the proof of this result, we have to consider the ``multi-field situation'' and we need to introduce some more notation first.

Let 
$$L\subseteq L_1',\ \ \ \ \ L\subseteq L_2',\ \ \ \ \ M\subseteq M_1',\ \ \ \ \ M\subseteq M_2'$$
be finite Galois extensions. We denote:
$$\widetilde{G}\left(\overline{L}'\right):=S_1\left(L_1'\right)\times S_2\left(L_2'\right)$$
and similarly for $\widetilde{H}\left(\overline{M}'\right)$. We further define:
$$\alpha_{\overline{L}'}:=\alpha_{L^{\alg}}|_{\widetilde{G}\left(\overline{L}'\right)}:\widetilde{G}\left(\overline{L}'\right)\to G\left(L^{\alg}\right),$$
$$\widetilde{G}_{\alpha}\left(\overline{L}'/L\right):=\left(\alpha_{\overline{L}'}\right)^{-1}(G(L))\leqslant \widetilde{G}\left(L\right),$$
$$G_{\overline{L}'}(L):=\alpha_{\overline{L}'}\left(\widetilde{G}_{\alpha}\left(\overline{L}'/L\right)\right)\leqslant G(L).$$
Note that all the groups and maps above are $L_2$-definable in the structure $(L,L)$, where $L_2$ is the two-sorted language of tuples of fields.

We similarly have $\beta_{\overline{M}'},\widetilde{H}_{\beta}\left(\overline{M}'/M\right),H_{\overline{M}'}(M)$ which are $L_2$-definable in the structure $(M,M)$.
\begin{prop}\label{mainmulti}
If the field maps $\varphi_1,\varphi_2:L\to M$ from Assumption \ref{assseveral} are elementary, then the induced group map $\varphi_G:G(L)\to G(M)$ is 1-elementary.
\end{prop}
\begin{proof}
Let $\phi_g(x)$ be a quantifier-free formula in the language of groups with the parameter $g\in G(L)$ and such that there is $h\in H(M)$ satisfying:
$$H(M)\models \phi_{\varphi_G(g)}(h)$$
(both $g$ and $h$ could be finite tuples for which the argument would be the same). There is a finite field extension $M\subseteq M'$ and $h',t\in \widetilde{H}(M')$ such that:
$$h=\beta_{M'}(h'),\ \ \ \ \ \varphi_G(g)=\beta_{M'}(t).$$
Since $h,\varphi_G(g)\in H_{M'}(M)$, $H(M)\models \phi_{\varphi_G(g)}(h)$ and $\phi$ is quantifier-free, we obtain that:
$$H_{M'}(M)\models \phi_{\varphi_G(g)}(h).$$
Using the fact that $H_{M'}(M)$ is $L_2$-definable in the structure $(M,M)$,
we get that there is an $L_2$-sentence over $(\varphi_1(L),\varphi_2(L))$ which is satisfied in $(M,M)$ and which says the following:
$$\text{``there are}\ (M_1/M,M_2/M)\ \text{finite Galois extensions such that}\ \ \ \ $$
$$\varphi_G(g)\in H_{M'}(M)\ \ \ \ \text{and}\ \ \ \ \ H_{M'}(M)\models (\exists x)\phi_{\varphi_G(g)}(x)\text{''}.$$
Since the field maps $\varphi_i:L\to M$ are elementary ($i=1,2$), the induced $L_2$-map:
$$(\varphi_1,\varphi_2):(L,L)\to (M,M)$$
is elementary as well (as in Remark \ref{ziegdisj}). Therefore, there are finite Galois extensions $L_1'/L,L_2'/L$ such that $g\in G_{\overline{L}'}(L)$ and
$$G_{\overline{L}'}(L)\models (\exists x)\phi_g(x).$$
Since $\phi$ is quantifier-free and $G_{\overline{L}'}(L)\leqslant G(L)$, we obtain that $G(L)\models (\exists x)\phi_g(x)$, which finishes the proof.
\end{proof}
The next result should be well-known but we could not find a direct reference, so we give a proof instead (and this proof fits to the context of this subsection).
\begin{fact}\label{coim}
If $L$ is infinite, then $\alpha_L(\widetilde{G}(L))=G(L)'$.
\end{fact}
\begin{proof}
By Theorem \ref{need}(2), we get $\widetilde{G}(L)=\widetilde{G}(L)'$, so
we obtain that $\alpha_L(\widetilde{G}(L))\subseteq G(L)'$.

For the other inclusion, we consider the group $\widetilde{G}_{\alpha}\left(L\right)$ from the proof of Proposition \ref{inducedg}.
We need the following.
\\
\\
{\bf Claim}
\\
$\widetilde{G}(L)$ is a normal subgroup of $\widetilde{G}_{\alpha}\left(L\right)$.
\begin{proof}[Proof of Claim]
By Claim 1 from the proof of Proposition \ref{inducedg}, it is enough to show that:
$$\widetilde{G}(L)\trianglelefteqslant \{g\in \widetilde{G}(L^{\alg})\ |\ \gal(L)g\subseteq N_{\alpha}(L^{\alg})g\}.$$

Take $x\in \widetilde{G}(L)$, $\sigma \in \gal(L)$ and $g\in \widetilde{G}(L^{\alg})$ such that $\gal(L)g\subseteq N_{\alpha}(L^{\alg})g$. We need to show that $\sigma(gxg^{-1})=gxg^{-1}$. By the assumption, there is $n\in N_{\alpha}(L^{\alg})$ such that $\sigma(g)=ng$ and we have:
\begin{IEEEeqnarray*}{rCl}
\sigma(gxg^{-1}) & = & \sigma(g)\sigma(x)\sigma(g)^{-1} \\
 & = & ngx(ng)^{-1}\\
 &= & gxg^{-1},
\end{IEEEeqnarray*}
 since $n\in Z(\widetilde{G}(L^{\alg}))$.
\end{proof}
By the Claim, $\alpha_L(\widetilde{G}(L))$ is normal in $G(L)'$, so (being an infinite subgroup of $G(L)'$ by Theorem \ref{zd}, since $L$ is infinite) coincides with $G(L)'$ by Theorem \ref{need}(1).
\end{proof}

\subsection{Automorphisms of finite commutative groups}\label{secaut}
This subsection may look a bit surprising, but the main result of it (Theorem \ref{subgroup}) will be crucial in the sequel (see Remark \ref{strproof}). For any commutative group $(W,+)$ and a natural number $k$, let us denote:
$$W[k]:=\{x\in W\ |\ kx=0\}.$$
\begin{theorem}\label{subgroup}
Let us fix the following:
\begin{itemize}
  \item $k,n\in \Nn$;
  \item finite cyclic groups $A_1,\ldots,A_n$;
  \item $A\leqslant A_1\times \ldots \times A_n$;
  \item $\alpha_i\in \mathrm{Aut}(A_i[k])$ for $i=1,\ldots ,n$
\end{itemize}
and assume that
$$\left(\alpha_1\times \ldots \times \alpha_n\right)\left(A[k]\right)=A[k].$$
Then, there are $\overline{\alpha}_i\in \mathrm{Aut}(A_i)$ extending $\alpha_i$ ($i=1,\ldots ,n$) and such that
$$\left(\overline{\alpha}_1\times \ldots \times \overline{\alpha}_n\right)\left(A\right)=A.$$
\end{theorem}
For the proof, we adopt our usual strategy, that is we assume that $n=2$ and leave an easy generalization to the case of an arbitrary finite $n$ to the reader. Let us assume then that $A\leqslant B\times C$. For some time we assume only that $B,C$ are commutative. Let
$$\pi_B:B\times C\to B,\ \ \ \ \pi_C:B\times C\to C$$
be the projection maps. We define:
$$B^A=\pi_B(A),\ \ C^A=\pi_C(A),\ \ B_A=A\cap \left(B\times \{0\}\right),\ \ C_A=A\cap \left(\{0\}\times C\right).$$
Then $A$ induces the following isomorphism:
$$f_A:B^A/B_A\to C^A/C_A.$$
\begin{lemma}\label{lfix1}
Assume that $\beta\in \mathrm{Aut}(B), \gamma\in \mathrm{Aut}(C)$. Then $(\beta\times \gamma)(A)=A$ if and only if
$$\beta(B_A)=B_A,\ \ \gamma(C^A)=C^A,\ \ \gamma(C_A)=C_A,\ \ (\beta\times \gamma)(\mathrm{graph}(f_A))=\mathrm{graph}(f_A).$$
\end{lemma}
\begin{proof}
The left-to-rigth implication is clear. For the other one, take $a\in A$. By the assumptions, we have:
$$a-(\beta\times \gamma)(a)\in B_A\times C_A.$$
Since $B_A\times C_A\subseteq A$, we obtain $(\beta\times \gamma)(a)\in A$.
\end{proof}
The next result is obvious.
\begin{lemma}\label{lfix2}
Assume that $\beta,\gamma,f\in \mathrm{Aut}(B)$. Then
$$(\beta\times \gamma)(\mathrm{graph}(f))=\mathrm{graph}(f)$$
 if and only if $f\circ \beta=\gamma\circ f$.

In particular, if $\mathrm{Aut}(B)$ is commutative then
$$(\beta\times \gamma)(\mathrm{graph}(f))=\mathrm{graph}(f)$$
if and only if $\beta=\gamma$.
\end{lemma}
Lemma \ref{lfix2} has the following specialization to the cyclic case.
\begin{lemma}\label{lfix3}
Assume that $B=B^A=\Zz/b\Zz$ and $C=C^A=\Zz/c\Zz$ for some $b,c\in \Nn$. Let
$$d:=|B/B_A|=|C/C_A|,$$
$\beta\in \mathrm{Aut}(B)$, and $\gamma\in \mathrm{Aut}(C)$. Then $(\beta\times \gamma)(A)=A$ if and only if
$$\beta(1+b\Zz) \equiv \gamma(1+c\Zz)\pmod{d}.$$
\end{lemma}
\begin{proof}
Since the groups $B$ and $C$ are cyclic, the conditions $\beta(B_A)=B_A,\gamma(C_A)=C_A$ are automatically fulfilled. Since homomorphisms are determined by their values on generators, we get the conclusion using Lemma \ref{lfix1} and Lemma \ref{lfix2}.
\end{proof}
We need now a very general result which must be folklore but we could not find a reference, so we give a proof instead. For a ring $T$, we denote by $T^*=\gm(T)$ its multiplicative group.
\begin{lemma}\label{ringepi}
Let $R,S$ be finite commutative rings with unity and $f:R\to S$ be a surjective ring homomorphism preserving unity. Then $f(R^*)=S^*$.
\end{lemma}
\begin{proof}
Since $R$ is finite, it is Artinian, so $R$ decomposes into a product of local Artinian rings (see \cite[Theorem 8.7.]{atmac}). Since any ideal in  a product of rings is a product of ideals, we can assume that $R$ is local, in which case it is easy to see that $f(R^*)=S^*$.
\end{proof}
\begin{remark}
We give some comments here about Lemma \ref{ringepi} and its proof.
\begin{enumerate}
\item Lemma \ref{ringepi} is not true for commutative rings in general, e.g. for the natural surjective homomorphism $\Zz\to \Zz/5\Zz$.

\item In the proof of Lemma \ref{ringepi}, only the Artinian assumption (rather than finiteness) was used and Lemma \ref{ringepi} is true for ``Artinian'' replacing ``finite''.

\item We are not aware if there is any more elementary argument for finite rings. But even for the rings of the form $\Zz/n\Zz$, it looks like the decomposition of $n$ into primes is crucial which exactly corresponds to decomposing an Artinian ring into a product of local Artinian rings.
\end{enumerate}
\end{remark}
\begin{proof}[Proof of Theorem \ref{subgroup}]
Let us take
$$B=\Zz/n\Zz,\ C=\Zz/m\Zz,\ k\in \Nn,\ \beta\in \mathrm{Aut}(B[k]),\ \gamma\in \mathrm{Aut}(C[k])$$
and assume that
$$\left(\beta\times \gamma\right)\left(A[k]\right)=A[k].$$
We can assume that $B=B^A,C=C^A$. Let
$$d:=|B/B_A|=|C/C_A|.$$
For any natural number $l$, we denote:
$$l_k:=\frac{l}{\mathrm{GCD}(l,k)}.$$
Then we have:
$$\left(\Zz/l\Zz\right)[k]\cong \Zz/l_k\Zz.$$
By Lemma \ref{lfix3}, we get:
$$\beta(1+n_k\Zz) \equiv \gamma(1+m_k\Zz)\pmod{d_k}.$$
Let us consider the following commutative diagram of finite commutative rings and the obvious surjective homomorphisms:
\begin{equation*}
 \xymatrix{ \Zz/n\Zz \ar[dd]^{} \ar[rrd]^{}  & &                       & & \Zz/m\Zz \ar[dd]^{} \ar[lld]^{} \\
 & & \Zz/d\Zz \ar[dd]^{}   & &  \\
\Zz/n_k\Zz  \ar[rrd]^{}  & &   & & \Zz/m_k\Zz  \ar[lld]^{} \\
 & & \Zz/d_k\Zz    & & }
\end{equation*}
Then we have:
$$(1+n_k\Zz,1+m_k\Zz)\in \left(\Zz/n_k\Zz\right)^*\times_{\left(\Zz/d_k\Zz\right)^*}\left(\Zz/m_k\Zz\right)^*$$
and therefore (using the interpretations from Lemma \ref{lfix3}), it is enough to show that the induced map:
$$\Phi^*:\left(\Zz/n\Zz\right)^*\times_{\left(\Zz/d\Zz\right)^*}\left(\Zz/m\Zz\right)^*\to  \left(\Zz/n_k\Zz\right)^*\times_{\left(\Zz/d_k\Zz\right)^*}\left(\Zz/m_k\Zz\right)^*$$
is onto. Using Lemma \ref{ringepi}, it is enough to show that the ring homomorphism
$$\Phi :\left(\Zz/n\Zz\right) \times_{\left(\Zz/d\Zz\right) }\left(\Zz/m\Zz\right) \to  \left(\Zz/n_k\Zz\right) \times_{\left(\Zz/d_k\Zz\right) }\left(\Zz/m_k\Zz\right)$$
is onto. To show it, we calculate the cardinalities of the domain, the range and the kernel of $\Phi$. Let us take $n',m',d'\in \Nn$ such that:
$$n=n'n_k,\ \ \ m=m'm_k,\ \ \ d=d'd_k.$$
Since we have the following additive isomorphism:
$$\ker(\Phi)\cong \left(\Zz/n'\Zz\right) \times_{\left(\Zz/d'\Zz\right) }\left(\Zz/m'\Zz\right),$$
we obtain:
$$|\mathrm{domain}(\Phi)|=\frac{nm}{d},\ \ \ |\mathrm{range}(\Phi)|=\frac{n_km_k}{d_k},\ \ \ |\ker(\Phi)|=\frac{n'm'}{d'}.$$
Putting it all together, we get:
$$|\mathrm{domain}(\Phi)|=|\mathrm{range}(\Phi)|\cdot |\ker(\Phi)|,$$
hence $\Phi$ must be onto.
\end{proof}

\section{The main result}
In this section, we show the main results of this paper. We start with the simple case, then we consider the Chevalley case and finally the semisimple case.

The next result is just a small part of the main theorem below (Theorem \ref{sscomm}), but we decided to include its proof separately, since it may be considered as a warm-up to the main proof and to the methods used there.
\begin{theorem}\label{smain}
Let $G$ be a split simply connected simple algebraic group over a field $K$ which is model complete. Then, the structure $(G(K),\cdot)$ is model complete (in the language of groups).
\end{theorem}
\begin{proof}
Let $H,N$ be models of $\mathrm{Th}(G(K),\cdot)$ and $f:H\to N$ be a monomorphism. We need to show that $f$ is elementary. Since any isomorphism is elementary and the composition of two elementary monomorphisms is again elementary, we will often replace $f$ with $h\circ f\circ g$, where $h$ and $g$ are group isomorphisms. By Remark \ref{mcspecial}, we can assume that $H$ and $N$ are special.

By Theorem \ref{will} applied to the interpretability functor
$$\mathrm{Th}_{\mathrm{Alg}_K}(K)\ni M\mapsto G(M)\in \mathrm{Groups},$$
there are field extensions $K\subseteq L$ and $L\subseteq M$ such that:
$$H\cong G(L),\ \ \ N\cong G(M),\ \ \ L\equiv K\equiv M \ \text{(in the language of $K$-algebras)}.$$
Therefore (by the ``we can often replace'' observation above), we can assume that $f:G(L)\to G(M)$.

By Theorem \ref{stein}, there is a field homomorphism $\varphi:L\to M$ and an isogeny (see the beginning of Section \ref{secclass} for the notation $G_L,G_M$)
$$\beta:  {}^{\varphi}(G_L)\to G_M$$
such that
$$f=\beta_{M}\circ \varphi_G.$$
Since $M$ is perfect (being model complete) we can use Fact \ref{coim} and Theorem \ref{need}(2), to get that $\beta_M$ is onto. Since $L\equiv M$, we obtain using Corollary \ref{ssdense} that
$$\varphi_G(Z(G(L)))=Z({}^{\varphi}(G_L)(M))\subseteq \ker(\beta_M).$$
Using that $\ker(f)$ is trivial, we obtain that $\beta_M$ is one-to-one, hence
$\beta_M$ is an isomorphism. Therefore, we can assume that $f=\varphi_G$, which is elementary, since the field $K$ is model complete and $L\equiv K\equiv M$.
\end{proof}
\begin{remark}
One could also see that $\beta_M$ is one-to-one using that $G_M$ is simply connected and then applying \cite[3.4]{Botits} to obtain that $\beta_M$ is an isomorphism.
\end{remark}
Below is the main result of this paper. The commutator group $G(K)'$ appearing there is the corresponding Chevalley group. These Chevalley groups are (abstractly) simple in the case when $G$ is adjoint. We give two separate proofs (for $G(K)$ and $G(K)'$), where the second proof uses the first one.
\begin{theorem}\label{sscomm}
Let $K$ be a model complete field and $G$ a semisimple and split algebraic group over $K$.
Then the groups $G(K)'$ and $G(K)$ are model complete.
\end{theorem}
\begin{proof}[Proof of Theorem \ref{sscomm} in the case of $G(K)'$]
Let $G_{\mathrm{sc}}$ be the ``universal cover'' of $G$ (see Notation \ref{scnot}). By Theorem \ref{prod}, we have the following decomposition:
$$G_{\mathrm{sc}}\cong S_1\times \ldots \times S_l,$$
where each $S_i$ is simple, simply connected and defined over $K$. Let
$$N:=\ker(\pi:G_{\mathrm{sc}}\to G)\leqslant Z(G_{\mathrm{sc}}).$$
The group $N(K)$ is finite and central in $G_{\mathrm{sc}}(K)$ and by Fact \ref{coim}, we obtain
$$G_{\mathrm{sc}}(K)/N(K)\cong G(K)'.$$
For simplicity of the presentation, we assume that $l=2$.

Take $G_1\equiv G(K)'\equiv G_2$ and let $f:G_1\to G_2$ be a monomorphism. As in the proof of Theorem \ref{smain} (note that Fact \ref{coim} implies that $G(K)'$ is definable), we can assume that $G_1=G(L)'$ and $G_2=G(M)'$ for $L\equiv K\equiv M$ (we consider here the theory of $K$ in the language of $K$-algebras).

Let $\pi:G_{\mathrm{sc}}(L)\to G(L)'$ denote the quotient map and
$$\widetilde{f}:=\pi_L\circ f:G_{\mathrm{sc}}(L)\to G(M)'.$$
By Theorem \ref{finalprod}, here are $\varphi_i:L\to M$ ($i=1,2$) and an isogeny
$$\beta:{}^{\varphi_1}S_1\times  {}^{\varphi_2}S_2 \to G$$
such that
$$\widetilde{f}=\beta_M\circ\left((\varphi_1)_{S_1}\times (\varphi_2)_{S_2}\right).$$
Let us introduce some notation corresponding to the one from Section \ref{secsf}:
$$\varphi_{\mathrm{sc}}:=(\varphi_1)_{S_1}\times (\varphi_2)_{S_2},\ \ \ {}^{\varphi}G_{\mathrm{sc}}:={}^{\varphi_1}S_1\times  {}^{\varphi_2}S_2.$$
We have the following commutative diagram:
\begin{equation*}
 \xymatrix{ G_{\mathrm{sc}}(L) \ar[d]^{\pi_L} \ar[rr]^{\varphi_{\mathrm{sc}}}  & & {}^{\varphi}G_{\mathrm{sc}}(M) \ar[d]^{\beta_M}
 \\ G(L) \ar[rr]^{f}  & &  G(M).}
\end{equation*}
{\bf Claim}
\\
We have:
$$\ker(\beta_M:{}^{\varphi}G_{\mathrm{sc}}(M)\to G(M))=\varphi_{\mathrm{sc}}(N(L)).$$
\begin{proof}[Proof of Claim]
Since $\pi_L\circ f=\beta_M\circ \varphi_{\mathrm{sc}}$, we obtain $\ker(\beta_M)\supseteq \varphi_{\mathrm{sc}}(N(L))$.

For the other inclusion, since the maps $\varphi_i$ are elementary, Fact \ref{rem: Hodges2} implies that the map $\varphi_{\mathrm{sc}}$ is elementary. As in the proof of Theorem \ref{smain}, this implies (since $Z(G_{\mathrm{sc}}(L))$ is finite) that:
$$\varphi_{\mathrm{sc}}\left(Z(G_{\mathrm{sc}}(L))\right)=Z({}^{\varphi}G_{\mathrm{sc}}(M)).$$
Hence we get:
$$\ker(\beta_M)\subseteq Z\left({}^{\varphi}G_{\mathrm{sc}}(M)\right)=\varphi_{\mathrm{sc}}\left(Z(G_{\mathrm{sc}}(L))\right)
\subseteq  \mathrm{im}(\varphi_{\mathrm{sc}}),$$
which gives the inclusion $\ker(\beta_M)\subseteq \varphi_{\mathrm{sc}}(N(L))$ using $\pi_L\circ f=\beta_M\circ \varphi_{\mathrm{sc}}$ (and the injectivity of $f$) again.
\end{proof}
Using Fact \ref{coim}, we get that the following homomorphisms
$$\pi_L:G_{\mathrm{sc}}(L)\to G(L)',\ \ \ \ \ \beta_M:{}^{\varphi}G_{\mathrm{sc}}(M)\to G(M)'$$
are onto. Therefore we have:
$$G(L)'\cong G_{\mathrm{sc}}(L)/N(L),\ \ \ \ G(M)'\cong {}^{\varphi}G_{\mathrm{sc}}(M)/\ker(\beta_M).$$
By Claim, we obtain:
$$f=\varphi_{\mathrm{sc}}/N(L):G(L)'\to G(M)'.$$
Since $\varphi_{\mathrm{sc}}$ is elementary (as in the proof of Claim) and $N(L)$ is definable (being finite), we get that $f$ is elementary.
\end{proof}
\begin{remark}\label{strproof}
Before the proof of Theorem \ref{sscomm} in the case of $G(K)$, we highlight its structure and extra difficulties. We have a control on the commutator group $G(K)'$ given by the group epimorphism $\pi_{\mathrm{sc}}:G_{\mathrm{sc}}(K)\to G(K)'$ using our understanding of the Chevalley case. In general, $G(K)$ need not be perfect, but $G(K^{\alg})$ is (see Theorem \ref{need}(4)). So, ``everything is fine'' over algebraically closed fields and (obviously) we can reach the algebraic closure of $K$ using finite extensions of $K$. Therefore, we need to be able to extend field homomorphisms in a very precise way, which is provided by Theorem \ref{subgroup}.
\end{remark}
\begin{proof}[Proof of Theorem \ref{sscomm} in the case of $G(K)$]
Let us take $G_1\equiv G(K)\equiv G_2$ and let $f:G_1\to G_2$ be a monomorphism. As above, we can assume that $G_1,G_2$ are special, so by Theorem \ref{will}, we can assume that $G_1=G(L)$ and $G_2=G(M)$ for some $K$-algebras $L,M$ such that $L\equiv K\equiv M$. As in the proof of the commutator case (applied to $f|_{G(L)'}$) we obtain field homomorphisms $\varphi_1,\varphi_2:L\to M$ and the commutative diagram:
\begin{equation*}
 \xymatrix{ G_{\mathrm{sc}}(L) \ar[d]^{\pi_L} \ar[rr]^{\varphi_{\mathrm{sc}}}  & & {}^{\varphi}G_{\mathrm{sc}}(M) \ar[d]^{\beta_M}
 \\ G(L) \ar[rr]^{f}  & &  G(M).}
\end{equation*}
such that
$$\varphi_{\mathrm{sc}}\left(\ker(\pi_L)\right)=\ker(\beta_M).$$
The next claim is crucial and in its proof we use the results from Subsection \ref{secaut}.
\\
\\
{\bf Claim 1}
\\
Each map $\varphi_i$ extends to a field homomorphism
$$\overline{\varphi}_i:L^{\alg}\to M^{\alg}$$
such that for the corresponding induced map
$$\overline{\varphi}_{\mathrm{sc}}:G\left(L^{\alg}\right)\to {}^{\varphi}G_{\mathrm{sc}}\left(M^{\alg}\right)$$
we have:
$$\overline{\varphi}_{\mathrm{sc}}\left(\ker(\pi_{L^{\alg}})\right)=\ker\left(\beta_{M^{\alg}}\right).$$
\begin{proof}[Proof of Claim 1]
We aim to use Theorem \ref{subgroup}. We need to perform several reductions first. Since
$$\ker(\pi)\leqslant Z\left(G_{\mathrm{sc}}\right),\ \ \ \ \ker(\beta)\leqslant Z\left({}^{\varphi}G_{\mathrm{sc}}\right),$$
we can replace $G_{\mathrm{sc}}$ with $Z(G_{\mathrm{sc}})$ and ${}^{\varphi}G_{\mathrm{sc}}$ with $Z\left({}^{\varphi}G_{\mathrm{sc}}\right)$.

By Theorem \ref{need}(5), for each simple algebraic group $S$, we have that $Z(S)\cong \mu_n$ for some $n>0$ or $Z(S)\cong \mu_2\times \mu_2$. Since for any field $F$, we have $\mu_2(F)=\mu_2(F^{\alg})$, we can assume that there are $n_i$ such that:
$$Z(S_i)\cong \mu_{n_i}\cong Z\left({}^{\varphi_i}S_i\right)$$
for $i=1,2$. Using the naturality of the above isomorphisms, for each
$$\overline{\varphi}_i:L^{\alg}\to M^{\alg}$$
as in the statement of Claim 1, we obtain the following commutative diagram:
\begin{equation*}
 \xymatrix{ Z\left(G_{\mathrm{sc}}\right)\left(L^{\alg}\right) \ar[d]^{\cong} \ar[rr]^{\overline{\varphi}_{\mathrm{sc}}}  & & Z\left({}^{\varphi}G_{\mathrm{sc}}\right)\left(M^{\alg}\right) \ar[d]^{\cong}
 \\ \left(\mu_{n_1}\times \mu_{n_2}\right)\left(L^{\alg}\right) \ar[rr]^{\overline{\varphi}}  & &  \left(\mu_{n_1}\times \mu_{n_2}\right)\left(M^{\alg}\right).}
\end{equation*}
Therefore, we can replace both of $Z\left(G_{\mathrm{sc}}\right)$ and $Z\left({}^{\varphi}G_{\mathrm{sc}}\right)$ with $\mu_{n_1}\times \mu_{n_2}$.

After all these replacements, we have $B\leqslant \left(\mu_{n_1}\times \mu_{n_2}\right)(K)$ corresponding both to $\ker(\pi_L)$ and $\ker(\beta_M)$, and $A\leqslant \left(\mu_{n_1}\times \mu_{n_2}\right)(K^{\alg})$ corresponding both to $\ker(\pi_{L^{\alg}})$ and $\ker(\beta_{M^{\alg}})$. We have also $\varphi(B)=B$ and our aim is to find the corresponding $\overline{\varphi}_i:L^{\alg}\to M^{\alg}$ ($i=1,2$) such that $\overline{\varphi}(A)=A$.

Let $F$ denote the prime subfield of $K$ and for any $n>0$ such that $\mathrm{char}(K)\nmid n$ let
$$F_n:=F(\zeta_n),$$
where $\zeta_n$ is the $n$-th primitive roof of unity in $F^{\alg}$. We fix $n$ as above such that:
$$\left(\mu_{n_1}\times \mu_{n_2}\right)(F_n)=\left(\mu_{n_1}\times \mu_{n_2}\right)\left(F^{\alg}\right)$$
and let $K_n:=K\cap F_n$. Then we have:
$$\left(\mu_{n_1}\times \mu_{n_2}\right)(K)=\left(\mu_{n_1}\times \mu_{n_2}\right)(K_n),$$
and there is $m>0$ such that (in the notation of Theorem \ref{subgroup}):
$$\left(\mu_{n_1}\times \mu_{n_2}\right)(K_n)=\left(\left(\mu_{n_1}\times \mu_{n_2}\right)(F_n)\right)[m],\ \ \ \ \ \ B=A[m].$$
For $i=1,2$ let us define
$$\alpha_i:=\varphi_i|_{\mu_{n_i}(K_n)}.$$
By Theorem \ref{subgroup}, there are $\overline{\alpha}_i\in \mathrm{Aut}(\mu_{n_i}(F_n))$ ($i=1,2$) such that
$$\left(\overline{\alpha}_1\times \overline{\alpha}_2\right)(A)=A.$$
It is enough to show now that each $\overline{\alpha}_i$ is a restriction of a field homomorphism $\overline{\varphi}_i:L^{\alg}\to M^{\alg}$ extending $\varphi$. Let $k_i>0$ be such that
$$\overline{\alpha}_i\left(\zeta_{n_i}\right)=\left(\zeta_{n_i}\right)^{k_i}.$$
Clearly, $k_i$ is relatively prime to $n_i$. Then, $\left(\zeta_{n_i}\right)^{k_i}$ is an $n_i$-th primitive root of unity as well, so there is $\beta_i\in \mathrm{Aut}(F_{n_i})$ such that $\beta_i(\zeta_{n_i})=\left(\zeta_{n_i}\right)^{k_i}$. Since the field extension $F\subseteq F_n$ is Abelian, there is $\gamma_i\in \mathrm{Aut}(F_{n})$ extending $\beta_i$ and (by the construction) extending $\varphi_i|_{K_n}$ as well.

Since $K\cap F_n=K_n$ and the extension $K_n\subseteq F_n$ is Galois, we obtain that $K$ and $F_n$ are linearly disjoint over $K_n$ and similarly both $L$ and $F_n$ are linearly disjoint over $K_n$, and $M$ and $F_n$ are linearly disjoint over $K_n$.
Therefore, we obtain:
$$LF_n\cong L\otimes_{K_n}F_n,\ \ \ \ \  MF_n\cong L\otimes_{K_n}M_n.$$
Hence, there is a field homomorphism $LF_n\to MF_n$ extending both $\varphi_i$ and $\gamma_i$ and it is enough to take any extension of this last field homomorphism to obtain our desired $\overline{\varphi}_i:L^{\alg}\to M^{\alg}$.
\end{proof}
By Claim 1, we are in the situation from Assumption \ref{assseveral}, so we obtain the corresponding group homomorphism $$\overline{\varphi}_G:G(L)\to G(M).$$
{\bf Claim 2}
\\
$\mathrm{im}(\overline{\varphi}_G)=\mathrm{im}(f)$.
\begin{proof}[Proof of Claim 2]
By Corollary \ref{ssdense} (since $|Z(G(L))|=|Z(G(M))|$ is finite) and using the set-up from Notation\ref{scnot}, there is the following quotient monomorphism:
$$f_{\mathrm{ad}}:\pi_{\mathrm{ad}}(G(L))\to \pi_{\mathrm{ad}}(G(M))$$
and its restriction to the commutator groups:
$$f_{\mathrm{ad}}|_{G_{\mathrm{ad}}(L)'}:G_{\mathrm{ad}}(L)'\to G_{\mathrm{ad}}(M)'.$$
Therefore, we obtain that
$$f_{\mathrm{ad}}|_{G_{\mathrm{ad}}(L)'}=\varphi_{G_{\mathrm{ad}}}|_{G_{\mathrm{ad}}(L)'}.$$
By Lemma \ref{coincide}, $f_{\mathrm{ad}}$ coincides with $\varphi_{G_{\mathrm{ad}}}$ on $\pi_{\mathrm{ad}}(G(L))$.
Therefore, for each $x\in G(L)$, we obtain:
$$\overline{\varphi}_G(x)f(x)^{-1}\in \ker\left(\pi_{\mathrm{ad}}\right)=Z(G(M)).$$
By Corollary \ref{ssdense}, $Z(G(M))\subseteq \mathrm{im}(f)$ and $Z(G(M))\subseteq \mathrm{im}(\overline{\varphi}_G)$, therefore we obtain that $\mathrm{im}(f)=\mathrm{im}(\overline{\varphi}_G)$.
\end{proof}
By Claim 2, Proposition \ref{mainmulti} and Lemma \ref{obv}, the map $f$ is 1-elementary which finishes the proof by Robinson's test.
\end{proof}
\begin{remark}
The following generalizations of Theorem \ref{sscomm} should not be difficult to prove using our methods.
\begin{enumerate}
  \item To other definable groups between $G(K)'$ and $G(K)$ for split semisimple algebraic groups $G$ defined over $K$.

  \item To semisimple algebraic groups $G$ which are not necessarily split but still satisfy the assumptions of the Borel-Tits Theorem as stated in \cite[Theorem 1.3]{Steinberg}. We may also make some assumptions like ``quasi-split'' to get the conditions of Theorem \ref{need}.
\end{enumerate}
However, if $G$ is simple algebraic and not split over $K$, then $G(K)$ can be even prosolvable (see \cite[Section 1.4.4]{plat93}). Therefore, the methods of this paper would not work in such a case and maybe one could find examples of a model complete field $K$ such that the group $G(K)$ in this case is \emph{not} model complete.

\end{remark}

\bibliographystyle{plain}
\bibliography{harvard}

\end{document}